\documentclass[10pt]{amsproc}
\usepackage{amsmath}
\usepackage{amsthm}
\usepackage{amssymb}
\usepackage{amsfonts}
\usepackage{tikz} 
\usepackage{yfonts}
\usepackage{mathtools}
\usepackage{cancel}
\usepackage{chemarrow}
\usetikzlibrary{matrix,arrows} 
\usepackage{amscd}
\newtheorem{theorem}{Theorem}[section]

\theoremstyle{definition}
\newtheorem{definition}[theorem]{Definition}

\theoremstyle{remark}
\newtheorem{remark}[theorem]{Remark}

\numberwithin{equation}{section}

\begin{document}

\title[Introduction to Higher rank Stable pairs and Virtual Localization]{An introduction to the theory of Higher rank stable pairs and Virtual localization}
\author{Artan Sheshmani}
\address{Max Plank Institute for Mathematics, Vivatsgasse 7, Bonn 53111}
\email{artan@mpim-bonn.mpg.de}
\urladdr{http://guests.mpim-bonn.mpg.de/artan/Welcome.html}
\date{June 3, 2012.}
\thanks{The author was supported in part by NSF Grants DMS 0244412, DMS 0555678 and DMS 08-38434 EMSW21-MCTP (R.E.G.S)}
\begin{abstract}
This article is an attempt to briefly introduce some of the results from \cite{a39} on development of a higher rank analog of the Pandharipande-Thomas theory  of stable pairs \cite{a17} on a Calabi-Yau threefold $X$.  More precisely,  we develop a moduli theory for highly frozen triples given by the data ${\mathcal O}_X^{\oplus r}(-n)\xrightarrow{\phi} F$ where $F$ is a sheaf of pure dimension $1$. The moduli space of such objects does not naturally determine an enumerative theory: that is, it does not naturally possess a perfect  
symmetric obstruction theory.  Instead, we build a zero-dimensional virtual fundamental class by hand, by truncating a deformation-obstruction theory coming from the moduli of objects in the derived  
category of $X$.  This yields the first deformation-theoretic construction of a higher-rank enumerative theory for Calabi-Yau threefolds.  We include the results of calculations in this enumerative theory for local ${\mathbb  
P}^1$ using the Graber-Pandharipande \cite{a25} virtual localization technique. We emphasize that in this article we merely include the statement of our theorems and illustrate the final result of some of the computations. The  proofs and more detailed calculations will appear elsewhere in a subsequent paper. 
\end{abstract}

\maketitle

\section*{Introduction}
The general philosophy behind the computation of Gromov-Witten invariants associated to a non-singular projective variety, $X$, is to do intersection theory over a moduli space that parametrizes curves embedded in $X$. In order to do intersection theory, one needs to compactify this moduli space. In 1992, Kontsevich \cite{a50} provided the moduli of stable maps as one example of such compactification. In this case, the resulting intersection numbers are rational; however,  there conjecturally exists a set of integers underlying these rational invariants which can be given an enumerative meaning. The second compactification was studied by Richard Thomas \cite{a20} which resulted in Donaldson-Thomas theory. The invariants computed in Donaldson-Thomas theory are integers but they do not directly count curves either, due to the fact that the one-dimensional subschemes of $X$ contain zero-dimensional subschemes which float around inside $X$. The third and the most recent compactification is due to Pandharipande and Thomas (PT) \cite{a17} which resulted in a new theory, known as the theory of stable pairs. Having fixed a nonsingular Calabi-Yau 3-fold $X$ with $\operatorname{H}^{1}(\mathcal{O}_{X})=0$, the objects in PT theory of stable pairs are given as tuples $(F,s)$, where $F$ is a coherent sheaf with one dimensional support given by a curve embedded inside $X$ and the section $s$ is defined by the map $s:\mathcal{O}_{X}\rightarrow F.$ Moreover, the stability condition for pairs means that the sheaf $F$ is pure with respect to its support and the cokernel of the morphism $s$ is given by a sheaf whose support is at most zero dimensional \cite{a17} (Lemma 1.3).
One would like to have a higher rank enumerative theory of Calabi-Yau 3-folds, compute new invariants and possibly relate them to physics (see for example \cite{a75} and \cite{a76}).  The present paper describes the step by step development of a higher rank analog of PT theory of stable pairs. We consider the higher rank objects in this theory called \textit{frozen triples} given as tuples $(E,F,\phi)$ where $E\cong \mathcal{O}_{X}^{\oplus r}(-n)$ (for some $r$) and $F$ is a pure sheaf with one dimensional support in $X$. Here $n\gg 0$ is a fixed large-enough integer which is introduced for technical deformation-theoretic reasons. We call the frozen triples of rank $r$ ``\textit{of type} $(P_{F},r)$" when $F$ has fixed Hilbert polynomial $P_{F}$. Note that given a frozen triple $(E,F,\phi)$ of type $(P_{F},r)$ we do not make a fixed choice of isomorphism $E\cong \mathcal{O}_{X}^{\oplus r}(-n)$. Having fixed such choice of isomorphism leads to a different moduli of closely related objects called \textit{highly frozen} triples given as tuples $(E,F,\phi,\psi)$ where $(E,F,\phi)$ have the same definition as before and $\psi:E\xrightarrow{\cong} \mathcal{O}_{X}^{\oplus r}(-n)$ is a fixed isomorphism.
\subsection*{Acknowledgement}
I would like to thank Sheldon Katz and Tom Nevins for their help and guidance. Their influence is all over my work. I would like to thank Richard Thomas for his continuous invaluable help and sharing with me his understanding of deformation obstruction theories in the derived category. Thanks to the organizers of the String-Math-2011 conference for providing me the opportunity to present this work and to have helpful discussions with the experts interested in this subject. 
\section{Moduli stacks}
The stability condition for frozen and highly frozen triples is compatible with PT stability of stable pairs \cite{a17} (Lemma 1.3). We call this stability condition $\tau'$-limit-stability or in short $\tau'$-stability. By definition \cite{a39} (Lemma 4.4), a   frozen (respectively highly frozen) triple  $(E,F,\phi)$ of type $(P_{F},r)$ is $\tau'$-limit-stable if and only if the map $E\xrightarrow{\phi}F$ has zero dimensional cokernel. We give a construction of the moduli space of $\tau'$-stable frozen and highly frozen triples as stacks. We show that these moduli stacks are given as algebraic stacks. More precisely, the moduli stack of $\tau'$-stable highly frozen triples of type $(P_{F},r)$ $\mathfrak{M}^{(P_{F},r,n)}_{s,\operatorname{HFT}}(\tau')$ is given as a Deligne-Mumford (DM) stack while the moduli stack of $\tau'$-stable frozen triples of type $(P_{F},r)$ $\mathfrak{M}^{(P_{F},r,n)}_{s,\operatorname{FT}}(\tau')$ is given as an Artin stack. We also show that $\mathfrak{M}^{(P_{F},r,n)}_{s,\operatorname{HFT}}(\tau')$ is a $\operatorname{GL}_{r}(\mathbb{C})$-torsor over $\mathfrak{M}^{(P_{F},r,n)}_{s,\operatorname{FT}}(\tau')$. We summarize these results as follows: 
\begin{theorem}\label{moduli-summary}
\cite{a39} (Proposition 5.5, Corollary 6.4, Theorems 6.2, 6.5)
\begin{enumerate}
\item Let $\mathfrak{S}_{s}^{(P_{F},r,n)}(\tau')$ be the stable locus of the parametrizing scheme of highly frozen triples of type $(P_{F},r)$ \cite{a39} (Section 6.1). Let $\left[\frac{\mathfrak{S}_{s}^{(P_{F},r,n)}(\tau')}{\operatorname{GL}(V_{F})}\right]$ be the stack-theoretic quotient of $\mathfrak{S}_{s}^{(P_{F},r,n)}(\tau')$ by $\operatorname{GL}(V_{F})$ where $V_{F}$ is defined as in \cite{a39} (Section 6.1). Then there exists an isomorphism of groupoids
\begin{equation} \label{hftisom}
\mathfrak{M}^{(P_{F},r,n)}_{s,\operatorname{HFT}}(\tau')\cong \left[\frac{\mathfrak{S}_{s}^{(P_{F},r,n)}(\tau')}{\operatorname{GL}(V_{F})}\right].
\end{equation}
\item Let $\left[\frac{\mathfrak{S}_{s}^{(P_{F},r,n)}(\tau')}{\operatorname{GL}_{r}(\mathbb{C})\times\operatorname{GL}(V_{F})}\right]$ be the stack-theoretic quotient of $\mathfrak{S}_{s}^{(P_{F},r,n)}(\tau')$ by $\operatorname{GL}_{r}(\mathbb{C})\times \operatorname{GL}(V_{F})$ where $V_{F}$. There exists an isomorphism of groupoids:$$\mathfrak{M}^{(P_{F},r,n)}_{s,\operatorname{FT}}(\tau')\cong \left[\frac{\mathfrak{S}_{s}^{(P_{F},r,n)}(\tau')}{\operatorname{GL}_{r}(\mathbb{C})\times \operatorname{GL}(V_{F})}\right]$$
\item The moduli stack $\mathfrak{M}^{(P_{F},r,n)}_{s,\operatorname{HFT}}(\tau')$ is a Deligne-Mumford (DM) stack
\item There exists a forgetful map $g':\mathfrak{M}^{(P_{F},r,n)}_{s,\operatorname{FT}}(\tau')\rightarrow \left[\frac{\operatorname{Spec}(\mathbb{C})}{\operatorname{GL}_{r}(\mathbb{C})}\right]$ such that the natural diagram:
\begin{equation}\label{diagram}
\begin{tikzpicture}
back line/.style={densely dotted}, 
cross line/.style={preaction={draw=white, -, 
line width=6pt}}] 
\matrix (m) [matrix of math nodes, 
row sep=2em, column sep=3.25em, 
text height=1.5ex, 
text depth=0.25ex]{ 
\mathfrak{M}^{(P_{F},r,n)}_{s,\operatorname{HFT}}(\tau')&pt=\operatorname{Spec}(\mathbb{C})\\
\mathfrak{M}^{(P_{F},r,n)}_{s,\operatorname{FT}}(\tau')&\mathcal{B}\operatorname{GL}_{r}(\mathbb{C})=\left[\frac{\operatorname{Spec}(\mathbb{C})}{\operatorname{GL}_{r}(\mathbb{C})}\right]\\};
\path[->]
(m-1-1) edge node [above] {$g$} (m-1-2)
(m-1-1) edge node [left] {$\pi$} (m-2-1)
(m-1-2) edge node [right] {$i$}(m-2-2)
(m-2-1) edge node [above] {$\acute{g}$} (m-2-2);
\end{tikzpicture},
\end{equation} 
is a fibered diagram in the category of stacks. In particular $\mathfrak{M}^{(P_{F},r,n)}_{s,\operatorname{HFT}}(\tau')$ is a $\operatorname{GL}_{r}(\mathbb{C})$-torsor over $\mathfrak{M}^{(P_{F},r,n)}_{s,\operatorname{FT}}(\tau')$. 
\end{enumerate}
\end{theorem}
\begin{remark}
The notion of $\tau'$-stability condition turns out to be a limiting GIT stability and thus one can apply the results of Wandel \cite{a62} (Section 3) to prove that the DM stack $\mathfrak{M}^{(P_{F},r,n)}_{s,\operatorname{HFT}}(\tau')$ has the stronger property of being given as a quasi-projective scheme. One can construct the higher rank theory of stable pairs over a nonsingular projective 3-fold or a noncompact 3-fold such as a toric variety.  However in the instance that the base 3-fold is chosen to be non-compact, in order to get well behaved moduli spaces, we require the condition that the one dimensional support of the sheaf $F$ appearing in a frozen or a highly frozen triple is compact. We will show later an example of this situation when $X$ is chosen to be a toric variety given by the total space of $\mathcal{O}_{\mathbb{P}^{1}}(-1)\oplus \mathcal{O}_{\mathbb{P}^{1}}(-1)\rightarrow \mathbb{P}^{1}$ (local $\mathbb{P}^{1}$).  
\end{remark}
\begin{remark}\label{classify}
If $X$ has a torus action, then the moduli stacks $\mathfrak{M}^{(P_{F},r,n)}_{s,\operatorname{HFT}}(\tau')$ and $\mathfrak{M}^{(P_{F},r,n)}_{s,\operatorname{FT}}(\tau')$ inherit it (once we have chosen an equivariant structure on $\mathcal{O}_{X}(-n)$) \cite{a39} (Section 12). Moreover, $\mathfrak{M}^{(P_{F},r,n)}_{s,\operatorname{HFT}}(\tau')$ has an extra ``\textit{non-geometric}" $\operatorname{T}_{0}:=(\mathbb{C}^{*})^{r}$ action by rescaling in components of $\mathcal{O}_{X}(-n)^{\oplus r}$ \cite{a39} (Section 13). Therefore, when $X$ is given as local $\mathbb{P}^{1}$, the natural ``\textit{geometric}" $\textbf{T}:=(\mathbb{C}^{*})^{3}$ action on $X$ induces an action on $\mathfrak{M}^{(P_{F},r,n)}_{s,\operatorname{FT}}(\tau')$ while $\mathfrak{M}^{(P_{F},r,n)}_{s,\operatorname{HFT}}(\tau')$ undergoes an action of $\textbf{G}:=(\mathbb{C}^{*})^{3}\times (\mathbb{C}^{*})^{r}$ which is due to the action of $\textbf{T}$ together with the additional non-geometric action of $\operatorname{T}_{0}$.
\end{remark}
\begin{remark}
As shown in \cite{a39} (Section 6), the construction of the moduli stacks of stable frozen and highly frozen triples depends on a choice of two fixed large enough integers $n\gg0$ and $n'\gg 0$. The first integer appears in the description of the triples $E\cong \mathcal{O}_{X}(-n)^{\oplus r}\rightarrow F$ and the second integer is the one for which $F(n')$ becomes globally generated and hence there exists a surjective map $V_{F}\otimes \mathcal{O}_{X}(-n')\twoheadrightarrow F$ for some $V_{F}$ (look at \cite{a39} (Section 6.1)). The existence of such $V_{F}$ is guaranteed according to Wandel \cite{a62} (Proposition 2.4) where he shows that given a bounded family of stable triples $E\rightarrow F$ there exists an integer $n'$ such that for every tuple $(E,F)$ appearing in the family $E(n')$ and $F(n')$ are globally generated over $X$. The fact that the sheaf $F(n')$ is globally generated for large enough values of $n'$ does not a priori imply that $\operatorname{H}^{i}(E_{2}(n))=0$ for all $i>0$ and our fixed choice of $n$. Hence we introduce the following definition:
\begin{definition}\label{open-sub}
\cite{a39} (Definition 7.1). Consider $\mathfrak{M}^{(P_{F},r,n)}_{s,\operatorname{HFT}}(\tau')$ and $\mathfrak{M}^{(P_{F},r,n)}_{s,\operatorname{FT}}(\tau')$. Now define the open substacks $\mathfrak{H}^{(P_{F},r,n)}_{s,\operatorname{HFT}}(\tau')\subset \mathfrak{M}^{(P_{F},r,n)}_{s,\operatorname{HFT}}(\tau')$ and $\mathfrak{H}^{(P_{F},r,n)}_{s,\operatorname{FT}}(\tau')\subset \mathfrak{M}^{(P_{F},r,n)}_{s,\operatorname{FT}}(\tau')$ as follows:
\begin{enumerate}
\item $\mathfrak{H}^{(P_{F},r,n)}_{s,\operatorname{HFT}}(\tau')=\{(E,F,\phi,\psi)\in \mathfrak{M}^{(P_{F},r,n)}_{s,\operatorname{HFT}}(\tau')\mid \operatorname{H}^{1}(F(n))=0\}$.
\item $\mathfrak{H}^{(P_{F},r,n)}_{s,\operatorname{FT}}(\tau')=\{(E,F,\phi)\in \mathfrak{M}^{(P_{F},r,n)}_{s,\operatorname{FT}}(\tau')\mid \operatorname{H}^{1}(F(n))=0\}$.
\end{enumerate}
\end{definition}
\end{remark}
From now on all our calculations are carried out over $\mathfrak{H}^{(P_{F},r,n)}_{s,\operatorname{HFT}}(\tau')$ and $\mathfrak{H}^{(P_{F},r,n)}_{s,\operatorname{FT}}(\tau')$ and the results in the following sections hold true for $\mathfrak{H}^{(P_{F},r,n)}_{s,\operatorname{HFT}}(\tau')$ and $\mathfrak{H}^{(P_{F},r,n)}_{s,\operatorname{FT}}(\tau')$ only. Also we assume that it is implicitly understood that in the following sections by the ``\textit{moduli stack of frozen or highly frozen triples}" we mean the open substack of the corresponding moduli stacks as defined in Definition \ref{open-sub}. As we will show later, though we prove our results over noncompact stacks, the torus fixed loci of these moduli stacks are compact which enable us to carry out localization computations over them.
\section{Deformation-obstruction theories}
For a 3-fold $X$ the natural deformation obstruction theories of stable frozen and highly frozen triples fail to provide well behaved complexes of correct amplitude over $\mathfrak{M}^{(P_{F},r,n)}_{s,\operatorname{FT}}(\tau')$ and $\mathfrak{M}^{(P_{F},r,n)}_{s,\operatorname{HFT}}(\tau')$ and they do not admit virtual cycles. We show that viewing the frozen and highly frozen triples as more complicated objects in $\mathcal{D}^{b}(X)$ given by $I^{\bullet}: E\rightarrow F$ and computing the fixed-determinant  obstruction theory of $I^{\bullet}$ will be the starting step in finding a well behaved deformation obstruction theory for the moduli stacks of frozen and highly frozen triples. This method has been successfully used in \cite{a17} to obtain an alternative candidate for the obstruction theory of the moduli of stable pairs. It is important to note that in the higher rank case, despite the fact that the object $I^{\bullet}$ (with the fixed determinant) in the derived category does not distinguish between a frozen or a highly frozen triple, its deformation space does. In other words, it can be shown that given a frozen triple $(E,F,\phi)$ and a highly frozen triple $(E,F,\phi,\psi)$, both associated to the same object $I^{\bullet}\in \mathcal{D}^{b}(X)$, the space of flat deformations of $(E,F,\phi)$ and $I^{\bullet}$ are equally governed by the group $\operatorname{Ext}^{1}(I^{\bullet},I^{\bullet})_{0}$ while the space of flat deformations of $(E,F,\phi,\psi)$ is not equal to that of $I^{\bullet}$. We summarize this remark as follows:
\begin{theorem}\cite{a39} (Propositions 7.2, 7.4 and Theorem 7.6)
\begin{enumerate}
\item Fix a map $f: S\rightarrow \mathfrak{H}^{(P_{F},r,n)}_{s,\operatorname{HFT}}(\tau')$. Let $S'$ be a square-zero extension of $S$ with ideal $\mathcal{I}$. Let $\mathcal{D}ef_{S}(S',\mathfrak{H}^{(P_{F},r,n)}_{s,\operatorname{HFT}}(\tau'))$ denote the deformation space of the map $f$ obtained by the set of possible deformations, $f':S'\rightarrow \mathfrak{H}^{(P_{F},r,n)}_{s,\operatorname{FT}}(\tau')$ such that $f'\mid_{S}=f$. The following statement is true:
\begin{equation}
\mathcal{D}ef_{S}(S',\mathfrak{H}^{(P_{F},r,n)}_{s,\operatorname{HFT}}(\tau'))\cong \operatorname{Hom}(I^{\bullet}_{S},F)\otimes \mathcal{I}
\end{equation}
\item Similarly for frozen triples let $f: S\rightarrow \mathfrak{H}^{(P_{F},r,n)}_{s,\operatorname{FT}}(\tau')$. Let $S'$ be a square-zero extension of $S$ with ideal $\mathcal{I}$. Let $\mathcal{D}ef_{S}(S',\mathfrak{H}^{(P_{F},r,n)}_{s,\operatorname{FT}}(\tau'))$ denote the deformation space of the map $f$ obtained by the set of possible deformations, $f':S'\rightarrow \mathfrak{H}^{(P_{F},r,n)}_{s,\operatorname{FT}}(\tau')$. The following statement is true:
\begin{align}
&
\mathcal{D}ef_{S}(S',\mathfrak{H}^{(P_{F},r,n)}_{s,\operatorname{FT}}(\tau'))\cong \operatorname{Ext}^{1}(I^{\bullet}_{S},I^{\bullet}_{S})_{0}\otimes \mathcal{I}.\notag\\
\end{align}
\end{enumerate}
\end{theorem}
We show that over $\mathfrak{H}^{(P_{F},r,n)}_{s,\operatorname{FT}}(\tau')$ deforming objects in the derived category leads to a 4-term deformation-obstruction complex of perfect amplitude $[-2,1]$:
\begin{theorem}\label{reldef-f}
\cite{a39} (Theorem 9.5). There exists a map in the derived category given by:$$R\pi_{\mathfrak{H}\ast}\left(R\mathcal{H}om(\mathbb{I}^{\bullet},\mathbb{I}^{\bullet})_{0}\otimes \pi_{X}^{*}\omega_{X}\right)[2]\xrightarrow{ob} \mathbb{L}^{\bullet}_{\mathfrak{H}^{(P_{F},r,n)}_{s,\operatorname{FT}}(\tau')}.$$ 
After suitable truncations, there exists a 4 term complex $\mathbb{E}^{\bullet}$ of locally free sheaves , such that $\mathbb{E}^{\bullet\vee}$ is self-symmetric of amplitude $[-2,1]$ and there exists a map in the derived category:
\begin{equation}\label{defobs-froz}
\mathbb{E}^{\bullet\vee}\xrightarrow{ob^{t}} \mathbb{L}^{\bullet}_{\mathfrak{H}^{(P_{F},r,n)}_{s,\operatorname{FT}}(\tau')},
\end{equation}
such that $h^{-1}(ob^{t})$ is surjective, and $h^{0}(ob^{t})$ and $h^{1}(ob^{t})$ are isomorphisms. Here $\mathbb{L}^{\bullet}_{\mathfrak{H}^{(P_{F},r,n)}_{s,\operatorname{FT}}(\tau')}$ stands for the truncated cotangent complex of the Artin stack $\mathfrak{H}^{(P_{F},r,n)}_{s,\operatorname{FT}}(\tau')$ which is of amplitude $[-1,1]$.
\end{theorem}
The computation of invariants over $\mathfrak{H}^{(P_{F},r,n)}_{s,\operatorname{FT}}(\tau')$ requires constructing a well-behaved virtual fundamental class.  Note that when $X$ is given as local $\mathbb{P}^{1}$, our method of computation is to use torus-equivariant cohomology and Graber Pandharipande (GP) virtual localization technique \cite{a25}. At the moment it is not clear to us how to classify the $\textbf{T}$-fixed loci of $\mathfrak{H}^{(P_{F},r,n)}_{s,\operatorname{FT}}(\tau')$. However this obstacle does not exist for the case of highly frozen triples. The existence of the additional non-geometric $T_{0}$ action (Remark \ref{classify}) over $\mathfrak{H}^{(P_{F},r,n)}_{s,\operatorname{HFT}}(\tau')$ makes it possible to classify the $\textbf{G}$-fixed locus of $\mathfrak{H}^{(P_{F},r,n)}_{s,\operatorname{FT}}(\tau')$ as a finite union of nonsingular compact components. Hence instead of developing a higher rank theory over $\mathfrak{H}^{(P_{F},r,n)}_{s,\operatorname{FT}}(\tau')$ we pull back the 4-term deformation obstruction complex in Theorem \ref{reldef-f} via the forgetful map $\pi$ in Theorem \ref{moduli-summary} (Diagram \eqref{diagram}) and try to construct a virtual fundamental class for $\mathfrak{H}^{(P_{F},r,n)}_{s,\operatorname{HFT}}(\tau')$. We summarize by saying that the construction of an enumerative theory over $\mathfrak{H}^{(P_{F},r,n)}_{s,\operatorname{HFT}}(\tau')$ has two advantages:
\begin{enumerate}
\item The construction of virtual fundamental classes and integration over DM stacks is in much more developed stage than over Artin stacks.
\item There exists a perfect classification of the torus fixed loci of highly frozen triples under the induced action of $\textbf{G}$ on $\mathfrak{H}^{(P_{F},r,n)}_{s,\operatorname{HFT}}(\tau')$ (Remark \ref{classify}) which makes it possible to do the computations using the virtual localization technique \cite{a25}.
\end{enumerate} 
Note that the complex $\pi^{*}\mathbb{E}^{\bullet\vee}$ is perfect of amplitude $[-2,1]$ and the main obstacle in constructing a well-behaved deformation obstruction theory over the DM stack $\mathfrak{H}^{(P_{F},r,n)}_{s,\operatorname{HFT}}(\tau')$ is to truncate $\pi^{*}\mathbb{E}^{\bullet\vee}$ into a 2-term complex and define (globally) a well-behaved deformation-obstruton theory of perfect amplitude $[-1,0]$. The simplest solution to this problem is to apply a cohomological truncation operation. Doing so requires obtaining a certain \textit{lifting} map from $g:\Omega_{\pi}\rightarrow \pi^{*}\mathbb{E}^{\bullet\vee}$ \cite{a39} (Proposition 9.12 and Lemma 9.13),  taking the mapping cone of this lift (and shifting by $-1$) and proving that the resulting complex satisfies the conditions of being a perfect obstruction theory for $\mathfrak{H}^{(P_{F},r,n)}_{s,\operatorname{HFT}}(\tau')$ \cite{a39} (Lemma 9.15). Here $\Omega_{\pi}$ is the relative cotangent sheaf of $\pi:\mathfrak{H}^{(P_{F},r,n)}_{s,\operatorname{HFT}}(\tau')\rightarrow \mathfrak{H}^{(P_{F},r,n)}_{s,\operatorname{FT}}(\tau')$. This procedure will remove the degree 1 term from the complex $\pi^{*}\mathbb{E}^{\bullet\vee}$. We also require to remove the degree $-2$ term of $\pi^{*}\mathbb{E}^{\bullet\vee}$ which is done by applying the same procedure to the map $g^{\vee}:\pi^{*}\mathbb{E}^{\bullet}\rightarrow \operatorname{T}_{\pi}$ \cite{a39} (Diagram 9.30) obtained from dualizing the map $g$. We finally obtain a local truncation of $\pi^{*}\mathbb{E}^{\bullet\vee}$ of perfect amplitude $[-1,0]$ which we denote by $\mathbb{G}^{\bullet}$. Assuming that $\pi^{*}\mathbb{E}^{\bullet\vee}$ is given by a 4 term complex of vector bundles:
\begin{equation}
\pi^{*}E^{-2}\rightarrow \pi^{*}E^{-1}\rightarrow \pi^{*}E^{0}\rightarrow \pi^{*}E^{1}
\end{equation}it can be seen from our construction \cite{a39} (Lemma 9.17) that locally the complex $\mathbb{G}^{\bullet}$ is given by$$\pi^{*}E^{-2}\xrightarrow{d'} \pi^{*}E^{-1}\oplus T_{\pi}\rightarrow \pi^{*}E^{0}\oplus \Omega_{\pi}\xrightarrow{d} \pi^{*}E^{1}$$which is quasi-isomorphic to a 2-term complex of vector bundles
\begin{equation}\label{G-dot}
\operatorname{Coker}(d')\rightarrow \operatorname{Ker}(d)
\end{equation}
concentrated in degree $-1$ and $0$.
\begin{remark}
In removing the degree -2 term of $\pi^{*}\mathbb{E}^{\bullet\vee}$ using the method described above, we exploit the fact that $\pi^{*}\mathbb{E}^{\bullet\vee}$ is a self-symmetric complex (look at \cite{a39} Diagram (9.30) for more discussion). 
\end{remark} 
The existence of the lifting map $g$ is guaranteed Zariski locally over $\mathfrak{H}^{(P_{F},r,n)}_{s,\operatorname{HFT}}(\tau')$ but not globally. Hence our strategy is to locally truncate $\pi^{*}\mathbb{E}^{\bullet\vee}$ as explained above, construct the corresponding local virtual cycles and glue the local cycles to define a globally-defined virtual fundamental class. Our main summarizing theorem of this section is as follows:
\begin{theorem}\label{2step-trunc}
\cite{a39} (Theorem 9.11). Consider the 4-term deformation obstruction theory $\mathbb{E}^{\bullet\vee}$ of perfect amplitude $[-2,1]$ over $\mathfrak{H}^{(P_{F},r,n)}_{s,\operatorname{FT}}(\tau')$.\\ 
1. Locally in the Zariski topology over $\mathfrak{H}^{(P_{F},r,n)}_{s,\operatorname{HFT}}(\tau')$ there exists a perfect two-term deformation obstruction theory of perfect amplitude $[-1,0]$ which is obtained from the suitable local truncation of the pullback $\pi^{*}\mathbb{E}^{\bullet\vee}$.\\ 
2. This local theory defines a globally well-behaved virtual fundamental class over $\mathfrak{H}^{(P_{F},r,n)}_{s,\operatorname{HFT}}(\tau')$.
\end{theorem}
Proving the second part of Theorem \ref{2step-trunc} requires an assumption which we explain next.
\section{Gluing the local virtual fundamental classes}
Let $\mathcal{U}=\coprod_{i}\mathcal{U}_{i}$ be an atlas of affine schemes for $\mathfrak{H}^{(P_{F},r,n)}_{s,\operatorname{HFT}}(\tau')$. Fix two open charts $\mathcal{U}_{\alpha}$ and $\mathcal{U}_{\beta}$ in $\mathcal{U}$. As we mentioned above, it is guaranteed that there exists lifting maps $g_{\alpha}:\Omega_{\pi}\mid_{\mathcal{U}_{\alpha}}\rightarrow \pi^{*}\mathbb{E}^{\bullet\vee}\mid_{\mathcal{U}_{\alpha}}$ and $g_{\beta}:\Omega_{\pi}\mid_{\mathcal{U}_{\beta}}\rightarrow \pi^{*}\mathbb{E}^{\bullet\vee}\mid_{\mathcal{U}_{\beta}}$ which lead to our desired locally defined deformation obstruction complexes over $\mathcal{U}_{\alpha}$ and $\mathcal{U}_{\beta}$ respectively. However, what essentially guarantees the existence of a globally well behaved virtual fundamental class is the compatibility of these two lifts over the intersection of $\mathcal{U}_{\alpha}$ and $\mathcal{U}_{\beta}$, i.e (roughly speaking) the gluing of $g_{\alpha}$ and $g_{\beta}$. Having fixed the lifting maps $g_{\alpha}$ and $g_{\beta}$ over $\mathcal{U}_{\alpha}$ and $\mathcal{U}_{\beta}$ enables one to apply the cohomological truncation operation as described in the previous section and obtain two locally truncated deformation obstruction theories of amplitude $[-1,0]$ which we denote by $\mathbb{G}^{\bullet}_{\alpha}$ and $\mathbb{G}^{\bullet}_{\beta}$. Now let $\mathcal{U}_{\alpha\beta}$ denote the intersection of $\mathcal{U}_{\alpha}$ and $\mathcal{U}_{\beta}$. Let $\mathbb{G}^{\bullet}_{\alpha}\mid_{\mathcal{U}_{\alpha\beta}}$ and $\mathbb{G}^{\bullet}_{\beta}\mid_{\mathcal{U}_{\alpha\beta}}$ denote the pullback of $\mathbb{G}^{\bullet}_{\alpha}$ and $\mathbb{G}^{\bullet}_{\beta}$ to $\mathcal{U}_{\alpha\beta}$. It is easy to see that via the restriction to $\mathcal{U}_{\alpha\beta}$ one obtains two maps in $\mathcal{D}^{b}(\mathcal{U}_{\alpha\beta})$ given by $\phi_{\alpha}\mid_{\mathcal{U}_{\alpha\beta}}:\mathbb{G}^{\bullet}_{\alpha}\mid_{\mathcal{U}_{\alpha\beta}}\rightarrow\mathbb{L}^{\bullet}_{\mathfrak{H}_{\operatorname{HFT}}}\mid_{\mathcal{U}_{\alpha\beta}}$ and $\phi_{\beta}\mid_{\mathcal{U}_{\alpha\beta}}:\mathbb{G}^{\bullet}_{\beta}\mid_{\mathcal{U}_{\alpha\beta}}\rightarrow\mathbb{L}^{\bullet}_{\mathfrak{H}_{\operatorname{HFT}}}\mid_{\mathcal{U}_{\alpha\beta}}$ . 
The gluability of the virtual cycles obtained from $\phi_{\alpha}$ and $\phi_{\beta}$ depends on whether $\phi_{\alpha}$ and $\phi_{\beta}$ satisfy the condition of ``\textit{giving rise to a semi-perfect obstruction theory}" over $\mathfrak{H}^{(P_{F},r,n)}_{s,\operatorname{HFT}}(\tau')$ in the sense of Liang-Li \cite{a70}:
 \begin{definition}
\cite{a70}(Definition 3.1). Let $X$ be a DM stack of finite type over an Artin stack $\mathcal{M}$. A semi perfect obstruction theory over $X\rightarrow \mathcal{M}$ consists 
of an \'{e}tale covering $\mathcal{U}=\coprod_{\alpha\in \Lambda}\mathcal{U}_{\alpha}$ of $X$ by schemes, and a truncated perfect relative obstruction theory $\phi_{\alpha}:\mathbb{G}^{\bullet}_{\alpha}\rightarrow \mathbb{L}^{\bullet}_{\mathcal{U}_{\alpha}/\mathcal{M}}$ for each $\alpha\in \Lambda$ such that
\begin{enumerate}
\item for each $\alpha,\beta$ in $\Lambda$ there is an isomorphism$$\psi_{\alpha\beta}:\operatorname{H}^{1}(\mathbb{G}^{\bullet}_{\alpha}\mid_{\mathcal{U}_{\alpha\beta}})\xrightarrow{\cong}\operatorname{H}^{1}( \mathbb{G}^{\bullet}_{\beta}\mid_{\mathcal{U}_{\alpha\beta}})$$so that the collection $(\operatorname{H}^{1}(\mathbb{G}^{\bullet}_{\alpha}),\psi_{\alpha\beta})$ forms a descent datum of sheaves.
\item For any pair  $\alpha,\beta\in\Lambda$ the obstruction theories $\phi_{\alpha}\mid_{\mathcal{U}_{\alpha\beta}}$ and $\phi_{\beta}\mid_{\mathcal{U}_{\alpha\beta}}$ are $\nu$-equivalent \cite{a70} (Definition 2.9).
\end{enumerate}
\end{definition}
It is easy to see that $\phi_{\alpha}$ and $\phi_{\beta}$ satisfy the condition (2) \cite{a39} (Proposition 10.4). However, (1) is a technical condition whose satisfaction depends on finding a specific choice of homotopical maps $h_{\alpha\beta}:g_{\alpha}\rightarrow g_{\beta}$ and $h^{\vee}_{\alpha\beta}:g^{\vee}_{\alpha}\rightarrow g^{\vee}_{\beta}$ where the maps $g$ and $g^{\vee}$ are the lifting maps defined in the previous section. In the present article and moreover in \cite{a39} we assume that such homotopical maps are given to us and we carry out our calculation of invariants over a specific toric Calabi-Yau 3-fold. However, finding rigorously the required homotopy maps needs further effort. We suspect that when $X$ is given as local $\mathbb{P}^{1}$ it  \textit{seems} possible to give a rigorous construction of such maps \cite{a77}.   
\begin{remark}\label{important}
It is important to note that essentially (as we will show later) the result of our calculations do not depend on the choice of the homotopy maps $h_{\alpha\beta}$ and $h^{\vee}_{\alpha\beta}$. In other words, the existence of the so-called well-defined homotopy maps will guarantee the existence of a theory of highly frozen triples for Calabi-Yau threefolds but no matter what choice of homotopy maps we make, it does not have any effect on the value of the numerical invariants which we calculate in this theory using the equivariant computations.    
\end{remark}
\section{Torus-fixed loci and virtual localization computations}
Let $X$ be given as the total space of $\mathcal{O}_{\mathbb{P}^{1}}(-1)\oplus \mathcal{O}_{\mathbb{P}^{1}}(-1)\rightarrow \mathbb{P}^{1}$ (local $\mathbb{P}^{1}$). Consider the ample line bundle over $X$ given by $\mathcal{O}_{X}(1)$. By our earlier notation, $\mathfrak{H}^{(P_{F},r,n)}_{s,\operatorname{HFT}}(\tau')$ and $\mathfrak{H}^{(P_{F},r,n)}_{s,\operatorname{FT}}(\tau')$ were defined as in Definition \ref{open-sub}.
It can be seen that when $X$ is given as local $\mathbb{P}^{1}$ then $\mathfrak{H}^{(P_{F},r,n)}_{s,\operatorname{HFT}}(\tau')=\mathfrak{M}^{(P_{F},r,n)}_{s,\operatorname{HFT}}(\tau')$ and $\mathfrak{H}^{(P_{F},r,n)}_{s,\operatorname{FT}}(\tau')=\mathfrak{M}^{(P_{F},r,n)}_{s,\operatorname{FT}}(\tau')$, hence in this section we switch back to our earlier notation and use $\mathfrak{M}^{(P_{F},r,n)}_{s,\operatorname{HFT}}(\tau')$ and $\mathfrak{M}^{(P_{F},r,n)}_{s,\operatorname{FT}}(\tau')$ instead. 
As we described (Remark \ref{classify}) there exists an induced action of $\textbf{G}:=\textbf{T}\times \operatorname{T}_{0}$ on $\mathfrak{M}^{(P_{F},r,n)}_{s,\operatorname{HFT}}(\tau')$. It can be shown that a torus fixed point in $\mathfrak{M}^{(P_{F},r,n)}_{s,\operatorname{HFT}}(\tau')$ corresponds to a $\textbf{G}$-equivariant highly frozen triple of type $(P_{F},r)$ \cite{a39} (propositions 12.2, 13.2). The key observation is that a $\textbf{G}$-equivariant highly frozen triple of rank $r$ is always written as a direct sum of $r$-copies of $\textbf{T}$-equivariant PT stable pairs \cite{a39} (Remark 13.3 and Remark 13.4):
\begin{equation}\label{split}
[\mathcal{O}^{\oplus r}_{X}(-n)\rightarrow F]^{\textbf{T}\times \operatorname{T}_{0}}\cong\bigoplus_{i=1}^{r}\left[\mathcal{O}_{X}(-n)\rightarrow F_{i}\right]^{\textbf{T}}.
\end{equation}
The consequence of this result is of significant importance since it enables one to immediately realize that the $\textbf{G}$-fixed loci of $\mathfrak{M}^{(P_{F},r,n)}_{s,\operatorname{HFT}}(\tau')$ are given as $r$-fold product of $\textbf{T}$-fixed loci of PT moduli space of stable pairs which are conjectured by Pandharipande and Thomas in \cite{a17} (Conjecture 2) to be nonsingular and compact. Hence, though our original moduli stack is constructed as a non-compact space, its torus-fixed locus is given as a finite union of compact and non-singular components. Let $\textbf{Q}$ denote a non-singular compact component of the torus fixed locus of $\mathfrak{M}^{(P_{F},r,n)}_{s,\operatorname{HFT}}(\tau')$. Let $\mathbb{G}^{\bullet\vee}_{\textbf{Q}}:G_{0,\textbf{Q}}\rightarrow G_{1,\textbf{Q}}$ be the dual of the restriction of the deformation obstruction complex in \eqref{G-dot} to $\textbf{Q}$. Using the methods described earlier, we construct the virtual fundamental class over all such $\textbf{Q}$ and obtain the virtual localization formula \cite{a39} (Equation 15.5):
\begin{equation}
\bigg[\mathfrak{M}^{r}_{\operatorname{s,HFT}}(\tau')\bigg]^{vir}=\sum_{\textbf{Q}\subset \mathfrak{M}^{r}_{\operatorname{s,HFT}}(\tau')}\iota_{\textbf{Q}*}\bigg(\frac{e(G_{1,\textbf{Q}})}{e(G_{0,\textbf{Q}})}\cdot e(T_{\textbf{Q}})\cap [\textbf{Q}]\bigg).
\end{equation}
Now we compute the difference $[G_{0,\textbf{Q}}]-[G_{1,\textbf{Q}}]$ in the $\textbf{G}$-equivariant $\mathcal{K}$-theory of $\textbf{Q}$. Consider a point $p\in \textbf{Q}$ represented by the complex $I^{\bullet\textbf{G}}:=[\mathcal{O}^{\oplus r}_{X}(-n)\rightarrow F]^{\textbf{G}}$. The difference $[G_{0,\textbf{Q}}]-[G_{1,\textbf{Q}}]$ over this point is the virtual tangent space at this point. We use the quasi isomorphism in  \eqref{G-dot} to compute the virtual tangent space:  
\begin{align}\label{virtg}
&
\mathcal{T}^{\textbf{Q}}_{I^{\bullet}}=[\operatorname{Coker}(d')]-[\operatorname{Ker}(d)]=\notag\\
&
\left([\pi^{*}E^{1}]-[\pi^{*}E^{0}]+[\pi^{*}E^{-1}]-[\pi^{*}E^{-2}]\right)+\left(\cancel{[T_{\pi}]}-\cancel{[\Omega_{\pi}]}\right),\notag\\
\end{align}
where $E^{i}$ for $i=-1,\cdots,2$ are the corresponding terms of $\mathbb{E}^{\bullet\vee}$ in Theorem \ref{reldef-f} and the cancellation in the second row is due to isomorphism of $\Omega_{\pi}$ and $T_{\pi}$ which is seen from their triviality \cite{a39} (Proposition 5.5 and exact triangle 9.21). By the construction of $\mathbb{E}^{\bullet\vee}$ and since the point $p\in \textbf{Q}$ is represented by the complex  $I^{\bullet\textbf{G}}$ the following identity holds true:
\begin{align}\label{virt-tang}
&
\mathcal{T}^{\textbf{Q}}_{I^{\bullet}}=\sum_{i=0}^{3}(-1)^{i+1}\cdot[\operatorname{Ext}^{i}(I^{\bullet\textbf{G}},I^{\bullet\textbf{G}})_{0}]=[\chi(\mathcal{O}_{X},\mathcal{O}_{X})]-[\chi(I^{\bullet\textbf{G}},I^{\bullet\textbf{G}})].\notag\\
\end{align}
The description of the virtual tangent space at a point plays an important role in the equivariant vertex and edge calculations \cite{a18} (Section 4.4 and Section 4.5). As it is seen, the right hand side of \eqref{virt-tang} is independent of choice of homotopy maps required for the existence of a globally well-behaved virtual fundamental class (Remark \ref{important}). Similar to computations in \cite{a18} (Section 4.3) we use \v{C}ech cohomology to compute the right hand side of \eqref{virt-tang} \cite{a39} (Section 15.2). Moreover, we give a computation of the equivariant vertex and edge for higher rank objects \cite{a39} (sections 15.2 and 15.3). 
\section{An Example} 
Assume that $X$ is given as local $\mathbb{P}^{1}$. There exists two affine patches $\mathcal{U}_{\alpha}$ and $\mathcal{U}_{\beta}$ covering $X$. The partitions associated to the Newton polyhedron of $X$ on each patch are given as three dimensional partitions with $\mu_{1}=(1),\mu_{2}=(0),\mu_{3}=(0)$ \cite{a18} (Example 4.9). We compute the vertex associated to the moduli stack of highly frozen triples of rank 2 given by $I^{\bullet}:=\mathcal{O}_{X}^{\oplus 2}(-n)\rightarrow F$.  Let $\mathcal{U}_{\alpha},\mathcal{U}_{\beta}$ denote affine open patches over the divisors $0,\infty$ on the base $\mathbb{P}^{1}$ respectively. Let $\mathbb{C}^{*}$ act on $\mathbb{C}^{4}$ by $t(x_{0},x_{1},x_{2},x_{3})=(tx_{0},tx_{1},t^{-1}x_{2},t^{-1}x_{3}).$ We identify $X$ as a quotient $X\cong (\mathbb{C}^{4}\backslash Z)\slash \mathbb{C}^{*}$ where $Z\subset \mathbb{C}^{4}$ is obtained by setting $x_{0}=x_{1}=0$. Let $([x_{0}:x_{1}],x_{2},x_{3})$ denote the coordiantes in $X$ where $[x_{0}:x_{1}]$ denote the homogeneous coordinates along the base $\mathbb{P}^{1}$ and $x_{2},x_{3}$ denote the fiber coordinates. Locally in the $\mathcal{U}_{\alpha}$ and $\mathcal{U}_{\beta}$ patches the defining coordinates are given as $(\frac{x_{1}}{x_{0}},x_{2}x_{0},x_{3}x_{0})$ and $(\frac{x_{0}}{x_{1}},x_{2}x_{1},x_{3}x_{1})$ respectively. Consider the $\mathcal{U}_{\alpha}$ patch. Let us denote the local coordinates in this patch by $(\tilde{x}_{1},\tilde{x}_{2},\tilde{x}_{3})$ where $\tilde{x}_{1}=\frac{x_{1}}{x_{0}},\tilde{x}_{2}=x_{2}x_{0},\tilde{x}_{3}=x_{3}x_{0}$. Now consider the action of $\textbf{T}=\mathbb{C}^{3}$ on $X$ where locally over $\mathcal{U}_{\alpha}$ is given by $(\lambda_{1},\lambda_{2},\lambda_{3})\cdot \tilde{x}_{i}=\lambda_{i}\cdot \tilde{x}_{i}.$ We identify an action of $(\mathbb{C}^{*})^{2}$ on $X$ which preserves the Calabi-Yau form by considering a subtorus $\operatorname{T}'\subset\textbf{T}$ such that
\begin{equation}\label{subtorus}
\operatorname{T}'=\{(\lambda_{1},\lambda_{2},\lambda_{3})\in \textbf{T}\mid \lambda_{1}\lambda_{2}\lambda_{3}=1\}.
\end{equation}
Let $\tilde{t}_{1},\cdots,\tilde{t}_{3}$ denote the characters corresponding to the action of $\lambda_{i}$. The equivariant vertex over the patch $\mathcal{U}_{\alpha}$ is obtained as follows:
\begin{align}\label{trace11}
& 
tr_{R-_{\chi((\mathbb{I}^{\bullet})_{\alpha},(\mathbb{I}^{\bullet})_{\alpha})}}=\textbf{F}^{\textbf{T}}_{\alpha}\cdot \frac{(w^{-1}_{1}+w^{-1}_{2})}{\tilde{t}^{n}_{1}}-\frac{\overline{\textbf{F}^{\textbf{T}}_{\alpha}}\cdot (w_{1}+w_{2})\cdot\tilde{t}^{n}_{1}}{\tilde{t}_{1}\tilde{t}_{2}\tilde{t}_{3}}\notag\\
&
+\textbf{F}^{\textbf{T}}_{\alpha}\overline{\textbf{F}^{\textbf{T}}_{\alpha}}\frac{(1-\tilde{t}_{1})(1-\tilde{t}_{2})(1-\tilde{t}_{3})}{\tilde{t}_{1}\tilde{t}_{2}\tilde{t}_{3}}+\frac{1-\frac{(w_{1}+w_{2})^{2}}{w_{1}w_{2}}}{(1-\tilde{t}_{1})(1-\tilde{t}_{2})(1-\tilde{t}_{3})},
\end{align}
where $(w_{1},w_{2})$ are defined as the weights of the action of $\operatorname{T}_{0}:=(\mathbb{C}^{*})^{\oplus 2}$ such that $w_{i}$ for $i=1,2$ are given by tuples $(1,0)$ and $(0,1)$. Here $\textbf{F}^{\textbf{T}}_{\alpha}$ denotes the equivariant character of the restriction (to the patch $\mathcal{U}_{\alpha}$) of the sheaf $F$. Similarly we compute $tr_{R-_{\chi((\mathbb{I}^{\bullet})_{\beta},(\mathbb{I}^{\bullet})_{\beta})}}$ \cite{a39} (Equation 16.8). The edge character is obtained as follows:
\begin{align}\label{trace22}
& 
tr_{R-_{\chi((\mathbb{I}^{\bullet})_{\alpha\beta},(\mathbb{I}^{\bullet})_{\alpha\beta})}}=\bigg(\textbf{F}^{\textbf{T}}_{\alpha\beta}\cdot (w^{-1}_{1}+w^{-1}_{2})-\frac{\overline{\textbf{F}^{\textbf{T}}_{\alpha\beta}}\cdot (w_{1}+w_{2})}{\tilde{t}_{2}\tilde{t}_{3}}+\notag\\
&
\textbf{F}^{\textbf{T}}_{\alpha\beta}\overline{\textbf{F}^{\textbf{T}}_{\alpha\beta}}\frac{((1-\tilde{t}_{2})(1-\tilde{t}_{3})}{\tilde{t}_{2}\tilde{t}_{3}}+\frac{1-\frac{(w_{1}+w_{2})^{2}}{w_{1}w_{2}}}{(1-\tilde{t}_{2})(1-\tilde{t}_{3})}\bigg)\delta(\tilde{t}_{1}).\notag\\
\end{align}
By compting the values of $\textbf{F}^{\textbf{T}}_{\alpha},\textbf{F}^{\textbf{T}}_{\beta},\textbf{F}^{\textbf{T}}_{\alpha\beta}$ we compute the vertex and edge characters in \eqref{trace11} and \eqref{trace22}. The $\textbf{G}$-character of the virtual tangent space in \eqref{virt-tang} is obtained by the following equation:
\begin{align}\label{virtan2}
&
tr_{R-\chi(I^{\bullet},I^{\bullet})}=tr_{R-\chi(\mathbb{I}^{\bullet}_{\alpha},\mathbb{I}^{\bullet}_{\alpha})}+ tr_{R-\chi(\mathbb{I}^{\bullet}_{\beta},\mathbb{I}^{\bullet}_{\beta})}-tr_{R-\chi(\mathbb{I}^{\bullet}_{\alpha\beta},\mathbb{I}^{\bullet}_{\alpha\beta})} 
\end{align}
Let $\textbf{Q}^{k}$ denote the $\textbf{G}$-fixed component of the moduli stack of rank 2 highly frozen triples over which the highly frozen triples $\mathcal{O}_{X}(-n)^{\oplus 2\textbf{G}}\xrightarrow{\phi} F^{\textbf{G}}$ satisfy the condition that $\operatorname{Length}(\operatorname{Coker}(\phi))=k$. By the definition of the equivariant vertex ($V_{\textbf{Q}}$) in \cite{a39} (Equation 15.34) the coefficient of the degree $k$ term in the vertex is obtained by the integral of the evaluation of the contribution of $V_{\textbf{Q}^{k}}$ on $\textbf{Q}^{k}$ \cite{a39} (Equation 16.12):
\begin{align}\label{lastnumber}
&
w(\textbf{Q}^{k})={\displaystyle{\int_{\textbf{Q}^{k}}}e(\operatorname{T}_{\textbf{Q}^{k}})e(-V_{\textbf{Q}^{k}})}=\notag\\
&
\displaystyle{\prod_{d_{1}+d_{2}=k}}\Bigg[\frac{\bigg(v_{1}(v_{2}-1)+\displaystyle{\prod_{i=0}^{d_{1}-1}}((i+n)s_{1})-(s_{2}+s_{3})\bigg)}{\bigg(v_{1}(1-v_{2})+\displaystyle{\prod_{i=1}^{d_{1}}}(-1)^{i}\cdot(i+n)s_{1}\bigg)}
\cdot\frac{\bigg(v_{2}(v_{1}-1)+\displaystyle{\prod_{i=0}^{d_{2}-1}}((i+n)s_{1})-(s_{2}+s_{3})\bigg)}{\bigg(v_{2}(1-v_{1})+\displaystyle{\prod_{i=1}^{d_{2}}}(-1)^{i}(i+n)s_{1}\bigg)}\Bigg]
\end{align}
where $\text{s}_{i}$ denote the equivariant characters corresponding to $\tilde{t}_{i}$. For similar discussions look at \cite{a18} (Section 4.7) as well as the calculation in \cite{a18} (Lemma 5). By the definition of the Calabi-Yau torus $\operatorname{T}'$ in \eqref{subtorus}, $s_{i}$ satisfy the property that $\text{s}_{1}+\text{s}_{2}+\text{s}_{3}=0$. Hence the generating series associated to the equivariant vertex of highly frozen triples of rank 2 is given by:
\begin{equation}
W^{\operatorname{HFT}}_{1,\emptyset,\emptyset}\mid_{r=2}=\left((1+q)^{\frac{(n+1)(s_{2}+s_{3})}{s_{1}}}\right)^{2}=\left(1+q\right)^{-2(n+1)}.
\end{equation} 
%This is an example of a special section head%
%%%%%%%%%%%%%%%%%%%%%%%%%%%%%%%%%%%%%%%%%%%%%%%%%%%%%%%%%%%%%%%%%%%%%%%%
%\footnote{Here is an example of a footnote. Notice that this footnote
%text is running on so that it can stand as an example of how a footnote
%with separate paragraphs should be written.
%\par
%And here is the beginning of the second paragraph.}%
%%%%%%%%%%%%%%%%%%%%%%%%%%%%%%%%%%%%%%%%%%%%%%%%%%%%%%%%%%%%%%%%%%%%%%%%
\bibliographystyle{amsalpha}
\bibliography{ref1}

\providecommand{\bysame}{\leavevmode\hbox to3em{\hrulefill}\thinspace}
\providecommand{\MR}{\relax\ifhmode\unskip\space\fi MR }
\providecommand{\MRhref}[2]{%
  \href{http://www.ams.org/mathscinet-getitem?mr=#1}{#2}
}
\providecommand{\href}[2]{#2}
\begin{thebibliography}{{Hua}11}

\bibitem[GP99]{a25}
T.~Graber and R.~Pandharipande, \emph{Localization of virtual classes}, Invent.
  Math. \textbf{135} (1999), 487--518.

\bibitem[{Hua}11]{a70}
{Huai-Liang Chang and Jun Li}, \emph{Semi-perfect obstruction theory and {DT}
  invariants of derived objects}, arXiv:1105.3261v1 (2011).

\bibitem[Kon92]{a50}
Maxim Kontsevich, \emph{Intersection theory on the moduli space of curves and
  the matrix {Airy} function}, Commun. Math. Phys. 147 (1992).

\bibitem[{Mic}09]{a75}
{Michele Cirafici, Annamaria Sinkovics, Richard J. Szabo}, \emph{Cohomological
  gauge theory, quiver matrix models and donaldson-thomas theory}, Nuclear
  physics B \textbf{809, Issue 3} (2009), 452--518.

\bibitem[PT09]{a18}
R.~Pandharipande and R.~P. Thomas, \emph{The 3-fold vertex via stable pairs},
  Geometry and Topology \textbf{13} (2009), 1835--1876.

\bibitem[{Ric}08]{a76}
{Michele Cirafici,}~{Annamaria Sinkovics,} {Richard J. Szabo},
  \emph{{Instantons and DonaldsonÐThomas invariants}}, Fortschritte der Physik
  \textbf{56, Issue 7-9} (2008), 849--855.

\bibitem[RR09]{a17}
{R. Pandharipande} and {R. P. Thomas}, \emph{{Curve} counting via stable pairs
  in the derived category}, Inventiones \textbf{178} (2009), 407--447.

\bibitem[She10]{a39}
Artan Sheshmani, \emph{Higher rank stable pairs and virtual localization},
  arXiv1011.6342 (2010).

\bibitem[{She}11]{a77}
{Artan} {Sheshmani}, \emph{{Semi}-perfect obstruction theories of higher rank
  stable pairs over local $\mathbb{P}^{1}$}, In preparation (2011).

\bibitem[Tho00]{a20}
R.~P. Thomas, \emph{A holomorphic {Casson} invariant for {Calabi}-{Yau}
  3-folds, and bundles on {K3} Þbrations}, J. Differential Geom. \textbf{54}
  (2000), 367--438.

\bibitem[Wan10]{a62}
Malte Wandel, \emph{Moduli spaces of stable pairs in {Donaldson}-{Thomas}
  theory}, arXiv:1011.3328v1 (2010).

\end{thebibliography}
%\begin{thebibliography}{A}

%\bibitem [A]{A} T. Aoki, \textit{Calcul exponentiel des op\'erateurs
%microdifferentiels d'ordre infini.} I, Ann. Inst. Fourier (Grenoble)
%\textbf{33} (1983), 227--250.

%\bibitem [B]{B} R. Brown, \textit{On a conjecture of Dirichlet},
%Amer. Math. Soc., Providence, RI, 1993.

%\bibitem [D]{D} R. A. DeVore, \textit{Approximation of functions},
%Proc. Sympos. Appl. Math., vol. 36,
%Amer. Math. Soc., Providence, RI, 1986, pp. 34--56.

%\end{thebibliography}

\end{document}